\documentclass[11pt,reqno]{amsart}
\usepackage{indentfirst, amssymb,amsmath,amsthm, mathrsfs}
\newtheorem*{theoA}{Theorem A}
\newtheorem*{theoB}{Theorem B}
\newtheorem*{theoC}{Theorem C}

\newtheorem{theorem}{Theorem}[section]
\newtheorem{lemma}{Lemma}[section]

\newtheorem{example}{Example}[section]

\newtheorem*{defiA}{Definition A}
\newtheorem*{defiB}{Definition B}
\newtheorem*{propA}{Proposition A}
\newtheorem*{propB}{Proposition B}
\newtheorem*{propC}{Proposition C}

\newcommand{\ol}{\overline}
\newcommand{\be}{\begin{equation}}
\newcommand{\ee}{\end{equation}}
\newcommand{\beas}{\begin{eqnarray*}}
\newcommand{\eeas}{\end{eqnarray*}}
\newcommand{\bea}{\begin{eqnarray}}
\newcommand{\eea}{\end{eqnarray}}
\newcommand{\bs}{\begin{small}}
\newcommand{\es}{\end{small}}
\renewcommand{\epsilon}{\varepsilon}
\numberwithin{equation}{section}
\newcommand{\lra}{\longrightarrow}
\begin{document}
\title[On meromorphic solutions of ...]
{On meromorphic solutions of certain Fermat-type difference and analogues equations concerning open problems}
\author[R. Mandal, R. Biswas and S. K. Guin]{Rajib Mandal, Raju Biswas and Sudip Kumar Guin}
\date{}
\address{Department of Mathematics, Raiganj University, Raiganj, West Bengal-733134, India.}
\email{rajibmathresearch@gmail.com}
\address{Department of Mathematics, Raiganj University, Raiganj, West Bengal-733134, India.}
\email{rajubiswasjanu02@gmail.com}
\address{Department of Mathematics, Raiganj University, Raiganj, West Bengal-733134, India.}
\email{sudipguin20@gmail.com}
\maketitle
\let\thefootnote\relax
\footnotetext{Type set by \AmS -\LaTeX}
\begin{abstract} In this paper, we have found that some certain Fermat-type shift and difference equations have the meromorphic solutions generated by Riccati type functions. 
Also we have solved the open problems posed by Liu and Yang (A note on meromorphic solutions of Fermat types equations, An. Stiint. Univ. Al. I. Cuza Lasi Mat. (N. S.), 62(2)(1), 317-325 (2016)). We have fortified the claims by some examples.\end{abstract}
\section{Introduction, Definitions and Results}
By a meromorphic (resp. entire) function, we shall always mean meromorphic (resp. entire) function over the complex plane $\mathbb{C}$.
Nevanlinna value distribution theory of meromorphic functions has been extensively applied to resolve growth (see e.g. \cite{7,8,18}) and solvability of meromorphic solutions
of linear and nonlinear differential equations (see e.g. \cite{40,41,42,17}). Let $f$ be a given meromorphic function on $\mathbb{C}$. We assume that the reader is familiar with the standard notations and results such as proximity function $m(r,f)$, counting function $N(r,f)$, characteristic function $T(r,f)$, the first and second main theorems, lemma on the logarithmic derivatives etc. of Nevanlinna theory (see e.g. \cite{7,8,18}). A meromorphic function $\alpha$ is said to be a small function of $f$, if $T(r,\alpha)=S(r,f)$, where $S(r,f)$ is used to denote any quantity that satisfies $S(r,f)=o(T(r,f))$ as $r\rightarrow \infty$, possibly outside of a set of $r$ of finite logarithmic measure. 
We denote by $\mathscr{S}(f)$ as the set of all small functions of $f$.\par
Let $a\in \mathscr{S}(f)\cap \mathscr{S}(g)$. For a meromorphic function $f$, if $f-a$ and $g-a$ have the same zeros with the same multiplicities, then we say that $f$ and $g$ share $a$ CM (counting multiplicities) and if we do not consider the multiplicities, then we say that $f$ and $g$ share $a$ IM (ignoring multiplicities).
We denote the order and the hyper order of a meromorphic function $f$ respectively by $\rho(f)$ and $\rho_2(f)$ such that 
\beas \rho(f)=\limsup\limits_{r\lra \infty}\frac{\log^+ T(r,f)}{\log r}\;\;\text{and}\;\;\rho_2(f)=\limsup\limits_{r\lra \infty}\frac{\log^+\log^+ T(r,f)}{\log r} .\eeas
Next, we recall the Hadamard's factorization theorem: Let $f(z)$ be meromorphic with $\rho(f )<+\infty$. Let $P_0(z)$ and $P_{\infty}(z)$ be the canonical products formed with the zeros and
poles of $f(z)$ in $\mathbb{C}\setminus\{0\}$ respectively. Let $c_mz^m$ with $c_m (\not=0)$ be the first non-vanishing term in the Laurent series of $f(z)$ near $0$. Then there exists a polynomial $Q(z)$ with
$\deg (Q)\leq \rho(f)$ such that $f(z)=z^me^{Q(z)}\frac{P_0(z)}{P_{\infty}(z)}$.\par
The next definitions are necessary for this paper.
\begin{defiA} Given a meromorphic function $f(z)$, $f(z+c)$ (resp. $f(qz+c)$) is called a shift (resp. $q$-shift) of $f$, where $c,q\in\mathbb{C}\setminus\{0\}$. Also for given a meromorphic function $f(z)$, $f(qz)$ is called a $q$-difference of $f$, where $q\in\mathbb{C}\setminus\{0\}$.\par
Given three meromorphic functions $f(z)$, $g(z)$ and $h(z)$, $f^n(z)+g^n(z)=h^n(z)$ is called a Fermat-type functional equation on $\mathbb{C}$, where $n\in\mathbb{N}$. 
Actually the functional equation is due to the assertion in \textrm{Fermat's Last Theorem} in $1637$ for the solutions of the Diophantine equation $x^n+y^n=z^n$ over some function fields, where $n\in\mathbb{N}$.\end{defiA}
\begin{defiB} $[$ Page-150, 228, \cite{kl}$]$
Let $f$ be a transcendental meromorphic function. Then functional equations of the forms
\bea\label{ekl} f(z+1)(\text{resp.}\;f(qz))=\frac{a_1(z)+b_1(z)f(z)}{c_1(z)+d_1(z)f(z)},\eea
where $a_1(z), b_1(z),c_1(z), d_1(z)\in \mathscr{S}(f)$ such that $a_1(z)d_1(z)-b_1(z)c_1(z)\not\equiv 0$, are called difference (resp. $q$-difference) Riccati equations.\\
In this paper, the equations similar to (\ref{ekl}), i.e., like $f(qz+c)=\frac{a_1(z)+b_1(z)f(z)}{c_1(z)+d_1(z)f(z)}$, we call as Riccati type equations, where $c,q(\not=0)\in\mathbb{C}$,  $a_1(z), b_1(z),c_1(z), d_1(z)\in \mathscr{S}(f)$ such that $a_1(z)d_1(z)-b_1(z)c_1(z)\not\equiv 0$.
\end{defiB}
We now consider the Fermat-type functional equation
\bea\label{eq1.1} f^n(z)+g^n(z)=1,\;\text{where}\;n\in\mathbb{N}.\eea 
The result due to Iyer \cite{2} is the gateway to find out the non-constant solutions of the Fermat-type functional equation (\ref{eq1.1}). \\
We summarize the classical results on the solutions of the functional equation (\ref{eq1.1}) in the following:
\begin{propA}
(i)\cite{22} The functional equation (\ref{eq1.1}) with $n=2$ has the non-constant entire solutions $f(z)=\cos(\eta(z))$ and $g(z)=\sin(\eta(z))$, where $\eta(z)$ is any entire function. No other solutions exist on $\mathbb{C}$.\\
(ii)\cite{3,22} For $n\geq 3$, there are no non-constant entire solutions of (\ref{eq1.1}) on $\mathbb{C}$.
\end{propA}
\begin{propB}
(i)\cite{22} The functional equation (\ref{eq1.1}) with $n=2$ has 
the non-constant meromorphic solutions 
$f=\frac{2\omega}{1+\omega^2}$ and $g=\frac{1-\omega^2}{1+\omega^2}$, where $\omega$ is an arbitrary meromorphic function on $\mathbb{C}$.\\ 
(ii)\cite{36,3} The functional equation (\ref{eq1.1}) with $n=3$ has
the non-constant meromorphic solutions 
$f=\frac{1}{2\wp(h)}\left(1+\frac{\wp'(h)}{\sqrt{3}}\right)$ and $g=\frac{\eta}{2\wp(h)}\left(1-\frac{\wp'(h)}{\sqrt{3}}\right)$, 
where $\eta$ is a cubic root of unity and $\wp(z)$ denotes the Weierstrass elliptic $\wp$-function with periods $\omega_1$ and $\omega_2$ is defined as
\beas \wp \left(z;\omega_1,\omega_2\right)=\frac{1}{z^2}+\sum\limits_{{\mu,\nu;\mu ^2+\nu ^2\not=0}} \left\{ \frac{1}{\left(z+\mu\omega_1+\nu\omega_2\right)^2}-\frac{1}{\left(\mu\omega_1+\nu\omega_2\right)^2}\right\},\eeas
which is even and satisfying, after appropriately choosing $\omega_1$ and $\omega_2$, $(\wp')^2=4\wp^3-1$.\\
(iii)\cite{3,22} For $n\geq 4$, there are no non-constant meromorphic solutions of (\ref{eq1.1}) on $\mathbb{C}$.
\end{propB}
In case \textrm{Proposition B(i)}, one may even obtain entire solutions of (\ref{eq1.1}), e.g., $f=\sin z$, $g=\cos z$, $\omega=\tan\frac{z}{2}$. In view of the transformation $\omega=\tan \left(\frac{h}{2}\right)$, where $h$ is an entire function, we see that in this case the functions $f=\frac{2\omega}{1+\omega^2}=\sin(h)$ and $g=\frac{1-\omega^2}{1+\omega^2}=\cos(h)$ are the only entire solutions.\\
We now focus on the non-constant solutions of the following functional equation 
\bea\label{equ1.1} a(z)f^n(z)+b(z)g^m(z)=1,\eea where $m$, $n\in\mathbb{N}$, $a(z)\in\mathscr{S}(f)$ and $b(z)\in\mathscr{S}(g)$.
In 1970, Yang \cite{16} investigated the functional equation (\ref{equ1.1}) as follows:
\begin{theoA} Let $m$, $n\in\mathbb{N}$ be satisfying $\frac{1}{m}+\frac{1}{n}<1$. Then there are no non-constant entire functions $f(z)$ and $g(z)$ satisfying (\ref{equ1.1}). 
\end{theoA}
Clearly the inequality in \textrm{Theorem A} holds for either $m\geq 2$, $n > 2$ or $m > 2$, $n\geq 2$. So, it is natural that the case $m=n=2$ can be treated when $f(z)$ and $g(z)$ have some special relationship in (\ref{equ1.1}), i.e., when $m=n=2$, the problem is still open.
This was the starting point of a new era about the solutions, mainly, entire solutions of functional equations like (\ref{equ1.1}).\par
Let $n\in\mathbb{N}$, $a$, $b_0$, $b_1$, $\cdots$, $b_{n-1}$ be polynomials, and $b_n\in\mathbb{C}\setminus\{0\}$. Let $L(f)=\sum\limits_{k=0}^nb_kf^{(k)}$ be a linear differential polynomial in $f$. In 2004, Yang and Li \cite{17} obtained that the solution of the Fermat-type equation
\bea\label{cp1}f^2+\left(L(f)\right)^2=a,\eea
must have the form $f(z)=\frac{1}{2}\left(P(z)e^{R(z)}+Q(z)e^{-R(z)}\right)$, where $P$, $Q$ and $R$ are polynomials with $PQ=a$.\par
Some meromorphic solutions of functional equation (\ref{equ1.1}) are found in \cite{14,34}.\par
{\bf Motivation:} The investigation to find out the solutions of the functional equation (\ref{equ1.1}) has gained a new dimension when $g(z)$ is replaced by the difference function $f(z+c)$ with finite order. 
The result due to Liu \cite{9} (mainly \textrm{Propositions 5.1} and \textrm{5.3}) is the gateway in this direction. Note that, there are many interesting contributions about the existences and forms of entire solutions with finite order of the equations similar to $f^2(z)+f^2(z+c)=1$. In this direction, all the finite order entire solutions are obtained by the use of \textrm{Hadamard's Factorization Theorem}. But as far our knowledge, no such method has been developed for finite order meromorphic solution. With this motivation, in this paper, we are actually interested to find out the finite order meromorphic solutions of three types of difference equations.\\
Next, we state the result on meromorphic solutions due to Liu et al. \cite{29} of the difference functional equation
\bea\label{e1}&& f^2(z)+f^2(z+c)=1,\;\text{where}\; c\in\mathbb{C}\setminus\{0\}\\\text{and}
\label{e1R}&&f^2(z)+f^2(qz)=1,\;\text{where}\; q(\not=0,1)\in\mathbb{C}\setminus\{0\}.\eea
\begin{propC} \cite{9} The meromorphic solutions of (\ref{e1}) (resp.  (\ref{e1R})) must satisfy $f(z)=\frac{1}{2}(h(z)+h^{-1}(z))$, where $h(z)$ is a meromorphic function satisfying one of the following two cases:
\begin{enumerate}
\item[(a)] $h(z+c)= - ih(z)$ $($resp. $h(qz)=-ih(z)$$)$;
\item[(b)] $h(z+c)h(z)=i$ $($resp. $h(qz)h(z)=i$$)$.\end{enumerate}
\end{propC}
\section{First main result}
In this paper, we consider the non-constant meromorphic functions with any order satisfying the Fermat-type functional equations (\ref{e1}), (\ref{e1R}) and 
\bea
\label{e1T}&&f^2(z)+f^2(qz+c)=1,\;\text{where}\; q(\not=1), c\in\mathbb{C}\setminus\{0\}\eea
 and show that the meromorphic solutions with any order of the equations can be generated by Riccati type meromorphic functions. Actually, we obtain the following result.
\begin{theorem}\label{th1} Let $f$ be a non-constant meromorphic function satisfying the Fermat-type functional equation (\ref{e1}) (resp. (\ref{e1R}) and (\ref{e1T})). Then $f(z)=\frac{2\omega(z)}{1+\omega^2(z)}$, where $\omega(z)$ is a non-constant meromorphic function satisfying the Riccati type equations either $\omega(z+c)(\text{resp.} \;\; \omega(qz)$ and $\omega(qz+c))=\frac{1+\omega(z)}{1-\omega(z)}$ or $\omega(z+c)(\text{resp.} \;\; \omega(qz)$ and $\omega(qz+c))=\frac{1-\omega(z)}{1+\omega(z)}$.\end{theorem}
The key tools in the proof of the results are \textrm{Propositions B} and \textrm{C}.
\subsection{Some examples of \textrm{Theorem \ref{th1}}}
 The following examples of \textrm{Theorem \ref{th1}} shows the importances of the results. Actually the results are the ways to construct the meromorphic solutions of the non-linear difference equations of types (\ref{e1}), (\ref{e1R}) and (\ref{e1T}).\\
We now exhibit the following examples for the Fermat type difference equation (\ref{e1}).
\begin{example} Let $\omega(z)=\frac{i\sqrt{-i}\cot (e^{4zi}+z)+1}{\sqrt{-i}\cot (e^{4zi}+z)+i}$ with $c=\frac{\pi}{2}$. Clearly then $\omega(z+c)=\frac{1-\omega(z)}{1+\omega(z)}$, which is a Riccati type equation and $f(z)=\frac{1}{2}\left(\sqrt{-i}\cot (e^{4zi}+z)+\frac{1}{\sqrt{-i}\cot (e^{4zi}+z)}\right)$ is a meromorphic solution of $f^2(z)+f^2\left(z+\frac{\pi}{2}\right)=1$. 
\end{example}
\begin{example} Let $\omega(z)=\frac{1+e^{i\left(e^z+\frac{z}{4i}\right)}}{i-ie^{i\left(e^z+\frac{z}{4i}\right)}}$ with $c=-2\pi i$. Note that $\omega(z)$ satisfies $\omega(z+c)=\frac{1+\omega(z)}{1-\omega(z)}$. We see that the function $f(z)=\frac{1}{2}\left(\frac{e^{i\left(e^z+\frac{z}{4i}\right)}}{i}+\frac{i}{e^{i\left(e^z+\frac{z}{4i}\right)}}\right)$ is a meromorphic solution of $f^2(z)+f^2\left(z-2\pi i\right)=1$. 
\end{example}
 For the Fermat type $q$-difference functional equation (\ref{e1R}), the following examples are necessary.
\begin{example} Let $\omega(z)=\frac{i(z+i)+e^{z^8}}{(z-i)+ie^{z^8}}$ with $q=-i$. It satisfies the corresponding Riccati type equation. As a result, the function $f(z)=\frac{1}{2}\left(\frac{z}{\left(e^{z^8}-1\right)}+\frac{\left(e^{z^8}-1\right)}{z}\right)$ is a meromorphic solution of $f^2(z)+f^2(-iz)=1$. 
\end{example}
\begin{example}Let $\omega(z)=\frac{i(1-\sqrt{i})\tan(\tan z)-i(1+\sqrt{i})}{(1+\sqrt{i})\tan(\tan z)-(1-\sqrt{i})}$ with $q=-1$ be satisfying the corresponding Riccati type equation. As a result, the meromorphic function $f(z)=\frac{1}{2}\left(\frac{\sqrt{i}(\tan(\tan z)-1)}{\tan(\tan z)+1}+\frac{\tan(\tan z)+1}{\sqrt{i}(\tan(\tan z)-1)}\right)$ is a solution of $f^2(z)+f^2(-z)=1$. 
\end{example}
Now we construct some examples for the Fermat type $q$-shift functional equation (\ref{e1T}) in the following.
\begin{example}Let $\omega(z)=\frac{i(1-\sqrt{i})\tan z-i(1+\sqrt{i})}{(1+\sqrt{i})\tan z-(1-\sqrt{i})}$ with $q=-1$ and $c=2\pi $. Note that $\omega(z)$ satisfies the corresponding Riccati type equation. As a result, the function  $f(z)=\frac{1}{2}\left(\frac{\sqrt{i}(\tan z-1)}{\tan z+1}+\frac{\tan z+1}{\sqrt{i}(\tan z-1)}\right)$ is a meromorphic solution of $f^2(z)+f^2(-z+2\pi )=1$. 
\end{example}
\begin{example} Let $\omega(z)=\frac{i(e^z-i)}{e^z+i}$ with $q=-1$ and $e^c=i$. It is easy to see that $\omega(z)$ satisfies $\omega(qz+c)=\frac{1-\omega(z)}{1+\omega(z)}$. Thus the function $f(z)=\frac{1}{2}\left(e^{-z}+e^z\right)$ is a meromorphic solution of $f^2(z)+f^2(-z+\frac{\pi}{2}i+2n\pi i)=1$, where $n \in \mathbb{Z}$. 
\end{example}
\subsection {Proof of Theorem \ref{th1}}  
\begin{proof}
The following three cases are for equations (\ref{e1}), (\ref{e1R}) and (\ref{e1T}) respectively.\\
{\bf Case 1.} We consider the Fermat type difference functional equation (\ref{e1}).
Note that, $f(z)\not\equiv f(z+c)$, otherwise, we get from (\ref{e1}) that $f$ is constant. Let the meromorphic solution of (\ref{e1}) be
\bea\label{e3b} f(z)=\frac{1}{2}\left(h(z)+\frac{1}{h(z)}\right),\eea
where $h(z)$ is a non-constant meromorphic function. 
In view of \textrm{Proposition B(i)}, we have
\bea\label{e3} &&f(z)=\frac{2\omega(z)}{1+\omega^2(z)}=\frac{1}{2}\left(\frac{i+\omega(z)}{i\omega(z)+1}+\frac{i-\omega(z)}{i\omega(z)-1}\right)\\
\text{and} \label{e4}&&f(z+c)=\frac{1-\omega^2(z)}{1+\omega^2(z)},\eea
where $\omega(z)$ is a non-constant meromorphic function.
From (\ref{e3b}) and (\ref{e3}), we set $h(z)=\frac{i+\omega(z)}{i\omega(z)+1}$.
Next in view of \textrm{Proposition C(i)}, we consider the following two cases.\\
{\bf Case 1.1.} When $h(z+c)=-ih(z)$. Now we easily deduce that 
\beas \frac{i+\omega(z+c)}{i\omega(z+c)+1}=-i\frac{i+\omega(z)}{i\omega(z)+1}
\Rightarrow \frac{-1+i\omega(z+c)}{i\omega(z+c)+1}=\frac{i+\omega(z)}{i\omega(z)+1}.\eeas
By componendo and dividendo rule, we easily get $\omega(z+c)=\frac{1+\omega(z)}{1-\omega(z)}$. Clearly this value satisfies (\ref{e4}), when it is calculated from (\ref{e3}).\\
{\bf Case 1.2.} When $h(z+c)h(z)=i$. Now we again easily deduce that 
\beas &&
\frac{i+\omega(z+c)}{i\omega(z+c)+1} \frac{i+\omega(z)}{i\omega(z)+1}=-\frac{1}{i}\nonumber\\\Rightarrow
&&\frac{1-i\omega(z+c)}{i\omega(z+c)+1}=\frac{i\omega(z)+1}{i+\omega(z)}.\eeas
By componendo and dividendo rule, we easily get $ \omega(z+c)=\frac{1-\omega(z)}{1+\omega(z)}$. Clearly this value satisfies (\ref{e4}), when it is again calculated from (\ref{e3}).\\
{\bf Case 2.} We consider the Fermat type $q$-difference functional equation (\ref{e1R}).
Using \textrm{Proposition C(ii)} instead of \textrm{Proposition C(i)} and proceeding as \textrm{Case 1}, we get the conclusions. \\
{\bf Case 3.} We consider the Fermat type $q$-shift functional equation (\ref{e1T}).
Note that, $f(z)\not\equiv f(qz+c)$, otherwise, we get from (\ref{e1T}) that $f$ is constant. Let the meromorphic solution of (\ref{e1T}) be
\bea\label{e3bT} f(z)=\frac{1}{2}\left(h(z)+\frac{1}{h(z)}\right)\;\text{and}\;f(qz+c)=\frac{1}{2i}\left(h(z)-\frac{1}{h(z)}\right),\eea
where $h(z)$ is a meromorphic function. Thus 
\beas&& \frac{1}{2i}\left(h(z)-\frac{1}{h(z)}\right)\equiv \frac{1}{2}\left(h(qz+c)+\frac{1}{h(qz+c)}\right)\\\Rightarrow
&&-ih(z)h(qz+c)[h(qz+c)+ih(z)]\equiv h(qz+c)+ih(z).\eeas
Thus, either $h(qz+c)\equiv -ih(z)$ or $h(z)h(qz+c)\equiv i$.
In view of \textrm{Proposition B(i)}, we know that
\bea\label{e3T} &&f(z)=\frac{2\omega(z)}{1+\omega^2(z)}=\frac{1}{2}\left(\frac{i+\omega(z)}{i\omega(z)+1}+\frac{i-\omega(z)}{i\omega(z)-1}\right)\\
\text{and} \label{e4T}&&f(qz+c)=\frac{1-\omega^2(z)}{1+\omega^2(z)},\eea
where $\omega(z)$ is a non-constant meromorphic function.
From (\ref{e3bT}) and (\ref{e3T}), we set $h(z)=\frac{i+\omega(z)}{i\omega(z)+1}$.
Next, we consider the following two cases.\\
{\bf Case 3.1.} When $h(qz+c)=-ih(z)$. Now we easily deduce that 
\beas &&\frac{i+\omega(qz+c)}{i\omega(qz+c)+1}=-i\frac{i+\omega(z)}{i\omega(z)+1}
\\\Rightarrow &&\frac{-1+i\omega(qz+c)}{i\omega(qz+c)+1}=\frac{i+\omega(z)}{i\omega(z)+1}.\eeas
By componendo and dividendo rule, we easily get $\omega(qz+c)=\frac{1+\omega(z)}{1-\omega(z)}$. Clearly this value satisfies (\ref{e4T}), when it is calculated from (\ref{e3T}).\\
{\bf Case 3.2.} When $h(qz+c)h(z)=i$. Now we again easily deduce that 
\beas &&
\frac{i+\omega(qz+c)}{i\omega(qz+c)+1} \frac{i+\omega(z)}{i\omega(z)+1}=-\frac{1}{i}\nonumber\\\Rightarrow
&&\frac{1-i\omega(qz+c)}{i\omega(qz+c)+1}=\frac{i\omega(z)+1}{i+\omega(z)}.\eeas
By componendo and dividendo rule, we easily get $ \omega(qz+c)=\frac{1-\omega(z)}{1+\omega(z)}$. Clearly this value satisfies (\ref{e4T}), when it is again calculated from (\ref{e3T}). 
This completes the proof.\end{proof}
\section{The Open Problems}
In the same paper, Liu and Yang \cite{29} considered another interesting Fermat-type functional equations, namely,
\bea\label{el1}&&f'(z)^2+f''(z)^2=f(z)^2\\\label{el2}\text{and}&&f\left(z+c_1\right)^2+f\left(z+c_2\right)^2=f(z)^2,\eea
where $c_1$ and $c_2$ are non-zero distinct constants. Obviously, the transcendental meromorphic solutions of (\ref{el1}) must be reduced to entire functions. Then they \cite{29} proved the following results.
\begin{theoB} The transcendental entire functions with finitely many zeros of (\ref{el1}) must be $f(z)=e^{az+b}$, where $a$ satisfies $a^4+a^2=1$.
\end{theoB}
\begin{theoC} The transcendental entire functions with finitely many zeros of (\ref{el2}) must be $f(z)=e^{a_1z+b_1}$, where $a_1$ satisfies $e^{2a_1c_1}+e^{2a_1c_2}=1$ and $b_1$ is a constant.\end{theoC}
Moreover, for further study for the solutions of (\ref{el1}) and (\ref{el2}), the authors \cite{29} raised the following questions as the open problems.\\
{\bf Open problem 1.} Whether there exists transcendental entire solutions of the functional equations (\ref{el1}) and (\ref{el2}) with infinitely many zeros ?\\
{\bf Open problem 2.} How to describe the transcendental meromorphic solutions of the functional equation (\ref{el2}) ?\\
As far we know, the above problems are unsolved till now.\\
For the entire solutions with infinitely many zeros of (\ref{el1}), we get the following result.
\begin{theorem}\label{th4} There is no transcendental entire function $f$ with infinitely many zeros satisfying (\ref{el1}).\end{theorem}
For the entire and meromorphic functions satisfying (\ref{el2}), we get the following result.
\begin{theorem}\label{th5} Let $f$ be a transcendental meromorphic (resp. entire) function satisfying (\ref{el2}). Then 
the following situations occur:
\item[(I)] If $\rho_2(f)<1$, then $f$ must be of finite order such that $f(z)=k_1e^{\frac{z}{c_1}\log M}\mathscr{F}(z)=k_2e^{\frac{z}{c_2}\log N}\mathscr{G}(z)$, where $\mathscr{F}$ and $\mathscr{G}$ are meromorphic periodic functions with periods $c_1$ and $c_2$ respectively and satisfying 
\beas \mathscr{F}\left(z+c_2\right)=\frac{k_2N}{k_1M^{\frac{c_2}{c_1}}}e^{\left(\frac{\log N}{c_2}-\frac{\log M}{c_1}\right)z}\mathscr{G}(z)\;\text{and}\; 
\mathscr{G}\left(z+c_1\right)=\frac{k_1M}{k_2N^{\frac{c_1}{c_2}}}e^{\left(\frac{\log M}{c_1}-\frac{\log N}{c_2}\right)z}\mathscr{F}(z),\eeas
where $M$, $N$ are distinct and $M^2+N^2=1$ with $k_1,k_2,M,N\in\mathbb{C}\setminus\{0\}$.
\item[(II)] If $\rho_2(f) \geq 1$ and, $\phi(z)$ and $\psi(z)$ are entire function and meromorphic function respectively such that $\psi(z)\equiv k_1\psi\left(z+c_1\right)\equiv k_2\psi\left(z+c_2\right)$, then $f$ must be of the form $f(z)=e^{\phi(z)}\psi(z)$, where $k_1,k_2\in\mathbb{C}\setminus\{0\}$ and, $e^{\phi(z)}$ and $\psi(z)$ are respectively the non-zero part and all-zero part of $f(z)$ with at least one of $e^{\phi(z)}$ and $\psi(z)$ is of hyper-order $\geq 1$. In particular, the following cases hold:
\item[(i)]  If $\phi(z)$ is both $c_1$ and $c_2$ periodic such that $\frac{c_1}{c_2}\in\mathbb{R}$, then $k_1^{-2}+k_2^{-2}=1$.
\item[(ii)] If $\phi(z)$ is neither $c_1$ nor $c_2$ periodic, then the following cases arise:
\begin{enumerate}
\item[(iia)] If $\phi(z)$ is a non-zero polynomial, then $\phi(z)=a_1z+a_0$, where $a_0,a_1\in\mathbb{C}$ such that $e^{2a_1c_1}+e^{2a_1c_2}=1$. Note that $\rho_2(\psi(z)) \geq 1$;
\item[(iib)] If $\phi(z)$ is a non-polynomial entire function, then $\phi(z)$ has one of the following forms:
$\phi(z)=\left\{\begin{array}{lll}
\Phi_1(z)+\frac{k}{c_1}z+c_3,&\text{when}\;\phi\left(z+c_1\right)-\phi(z)\equiv k,\\
\Phi_2(z)+\frac{k}{c_2}z+c_4,&\text{when}\;\phi\left(z+c_2\right)-\phi(z)\equiv k,\\
\Phi_3(z)+\frac{k}{c_1-c_2}z+c_5,&\text{when}\;\phi\left(z+c_1\right)-\phi(z+c_2)\equiv k,\end{array}\right.$\\
where $\Phi_1(z),\Phi_2(z)$ and $\Phi_3(z)$ are transcendental entire periodic functions with periods $c_1,c_2,c_1-c_2$ respectively and $k(\not=0),c_3,c_4,c_5\in\mathbb{C}$.
\end{enumerate}
\item[(iii)] If $\phi(z)$ is either $c_1$ or $c_2$ periodic, then either $\phi\left(z+c_2\right)-\phi(z)\equiv k$ with $k_1^{-2}\not=1$ or $\phi\left(z+c_1\right)-\phi(z)\equiv k$ with $k_1^{-2}\not=1$, where $k(\not=0)\in\mathbb{C}$ and the similar conclusions as (iib) hold.
\end{theorem}
Clearly we have solved the {\bf Open problems 1} and {\bf 2} in \textrm{Theorems \ref{th4}} and \textrm{\ref{th5}}.\\
The key tools in the proof of the results are \textrm{Hadamard's Factorization Theorem} and the core part of Nevanlinna's theory.
\subsection{Some examples of \textrm{Theorem \ref{th5}}}
Actually, the entire and meromorphic functions with infinitely many zeros satisfying (\ref{el2}) exist. We exhibit the following examples for the claim.
\begin{example}
Let $k_1=k_2$, $\mathscr{F}(z)=\mathscr{G}(z)=\tan{\pi z}$, $c_1=2$, $c_2=4$, $M={\left(\frac{\sqrt{5}-1}{2}\right)}^{\frac{1}{2}}$ and $N=\frac{\sqrt{5}-1}{2}$ in \textrm{Theorem \ref{th5}}. Then we see that $f(z)=k_1e^{\frac{z}{c_1}\log M}\tan{\pi z}$ is a finite order transcendental meromorphic solution of (\ref{el2}) and satisfying the conditions in \textrm{Theorem \ref{th5}}.
\end{example}
\begin{example}
Let $a_1=\frac{1}{4}\log{\frac{-\sqrt{5}-1}{2}}$,  $c_1=2$, $c_2=4$ and $\psi(z)=\cot(\pi z+k_0)$. Clearly then $f(z)=e^{\frac{1}{4}\log{\frac{-\sqrt{5}-1}{2}}}\cot(\pi z+k)$ is a transcendental meromorphic solution of (\ref{el2}) and $e^{2a_1c_1}+e^{2a_1c_2}=1$, where $a_0,k_0\in\mathbb{C}$.
\end{example}
\begin{example}
Let $f(z)=e^{\cosh {2z}+\frac{k}{c_2}z+c_4}\psi(z)$ and $k=\log{\frac{\sqrt{5}-1}{2}}$, where $c_3\in\mathbb{C}$ and $\psi(z)=\sinh {2z}$ with periods $c_1=\pi i$ and $c_2=2\pi i$. We see that $f$ is a transcendental entire solution of (\ref{el2}) with infinitely many zeros and $e^{k}+e^{2k}=1$.
\end{example}
\begin{example}
Let $k=\frac{\pi}{6}i$ and $\psi(z)=\sin z$ with periods $c_1=\pi$ and $c_2=-\pi$. Then $f(z)=e^{\cos{(\sin z)}+\frac{k}{c_2}z+c_4}\psi(z)$ is a transcendental entire solution of (\ref{el2}) with infinitely many zeros and $e^{-2k}+e^{2k}=1$, where $c_4\in\mathbb{C}$.
\end{example}

\subsection {The technical lemmas} 
The following lemmas are used to prove the {\bf Open problems 1} and {\bf 2}. 
\begin{lemma}\label{lem9}\cite{18} Suppose that $f_j(z)$ are meromorphic functions and $g_k$ are entire functions ($1\leq j\leq  n$, $n\geq 2$) satisfying the following conditions:
\begin{enumerate}
\item[(i)] $\sum\limits_{j=1}^{n}f_j(z)e^{g_j(z)}\equiv 0$,
\item[(ii)] $g_j(z)-g_k(z)$ are not constants for $1\leq j< k\leq  n$,
\item[(iii)] for $1\leq j\leq  n$, $1\leq h<k\leq  n$, $T\left(r,f_j\right)=o\{T\left(r,e^{g_h-g_k}\right)\}$ ($r\to \infty$, $r\not\in E$).
\end{enumerate}
Then $f_j(z)\equiv 0$ for $1\leq j\leq n$.\end{lemma}
\begin{lemma}\label{lem6}\cite{18} Let $f_i\in\mathscr{M}$ for $1\leq i \leq n, n\geq 3$ and are not constants except for $f_n$. Also, let $\sum_{i=1}^nf_i\equiv 1$. If $f_n\not\equiv 0$ and 
\beas \sum_{i=1}^n N\left(r,0; f_i\right)+(n-1)\sum_{i=1}^n\ol N\left(r, f_i\right)<(\lambda+o(1))T\left(r, f_k\right),\eeas
where $\lambda<1$ and $1\leq k\leq n-1$. Then $f_n\equiv 1$.\end{lemma}
The following basic inequalities, by $[$\cite{50e}, \textrm{Lemma} 8.3$]$, are frequently used in value distribution theory for differences.
\begin{lemma}\label{lem8l}\cite{50e}
Let $f(z)$ be a non-constant meromorphic function with hyper-order less than $1$, $c \in\mathbb{C}$. Then, \beas\label{ep1m} &&N(r,0;f(z+c))\leq N(r,0;f(z))+S(r,f),\\&&N(r,f(z+c))\leq N(r,f)+S(r,f),\nonumber\\ 
\label{ep3} &&\ol N(r,0;f(z+c))\leq \ol N(r,0;f(z))+S(r,f)\\\text{and}&&\ol N(r,f(z+c))\leq \ol N(r,f)+S(r,f).\eeas \end{lemma}
\subsection {Proofs of the Open problems 1 and 2}  
\begin{proof}[\bf{Proof of Theorem \ref{th4}}] Let $f(z)$ be a transcendental entire function with infinitely many zeros satisfying (\ref{el1}). 
Rewritting (\ref{el1}), we get
\bea\label{el4s} \left(\frac{f'\left(z\right)}{f(z)}\right)^2+\left(\frac{f''\left(z\right)}{f(z)}\right)^2=1.\eea 
Clearly both of $\frac{f'\left(z\right)}{f(z)}$ and $\frac{f''\left(z\right)}{f(z)}$ are either entire functions or non-entire meromorphic functions. 
Now we consider the following cases.\\
{\bf Case 1.} Let $\frac{f'\left(z\right)}{f(z)}$ and $\frac{f''\left(z\right)}{f(z)}$ be both entire functions. In view of \textrm{Proposition A (i)}, we claim that
\bea\label{el7as} \frac{f'\left(z\right)}{f(z)}=\cos(\eta(z))\;\;
\text{and}\;\frac{f''\left(z\right)}{f(z)}=\sin(\eta(z))
,\eea
where $\eta(z)$ is an entire function on $\mathbb{C}$. If we integrate the first part in (\ref{el7as}), then we see that $f(z)$ has no zeros and this is not possible.\\
{\bf Case 2.} Let $\frac{f'\left(z\right)}{f(z)}$ and $\frac{f''\left(z\right)}{f(z)}$ be both non-entire meromorphic functions. In view of \textrm{Proposition B (i)}, we claim that
\bea\label{el7es} \frac{f'\left(z\right)}{f(z)}=\frac{2\omega(z)}{1+\omega^2(z)}\;\;
\text{and}\;\;\frac{f''\left(z\right)}{f(z)}=\frac{1-\omega^2(z)}{1+\omega^2(z)}
,\eea
where $\omega(z)$ is a non-constant meromorphic function on $\mathbb{C}$. Now we consider the following cases.\\
{\bf Sub-case 2.1} Let $\int \frac{2\omega(z)}{1+\omega^2(z)}$ be not of the form $\log \psi(z)$ for some transcendental entire function $\psi(z)$. If we integrate the first part in (\ref{el7es}), then we see that $f(z)$ has no zeros and this is not possible.\\
{\bf Sub-case 2.2} Let $\int \frac{2\omega(z)}{1+\omega^2(z)}$ be of the form $\log \psi(z)$ for some transcendental entire function $\psi(z)$ with infinitely many zeros. From (\ref{el7es}), we easily get
\beas &&f''(z)=2\frac{\left(1+\omega^2(z)\right)\omega'(z)-2\omega^2(z)\omega'(z)}{\left(1+\omega^2(z)\right)^2}f(z)+\frac{2\omega(z)}{1+\omega^2(z)}f'(z)\\
\Rightarrow&&f''(z)=2\frac{1-\omega^2(z)}{\left(1+\omega^2(z)\right)^2}\omega'(z)f(z)+\left(\frac{2\omega(z)}{1+\omega^2(z)}\right)^2f(z)\\
\Rightarrow&&\frac{f''(z)}{f(z)}=2\frac{1-\omega^2(z)}{\left(1+\omega^2(z)\right)^2}\omega'(z)+\left(\frac{2\omega(z)}{1+\omega^2(z)}\right)^2\\
\Rightarrow&&\frac{1-\omega^2(z)}{1+\omega^2(z)}\equiv 2\frac{1-\omega^2(z)}{\left(1+\omega^2(z)\right)^2}\omega'(z)+\left(\frac{2\omega(z)}{1+\omega^2(z)}\right)^2\\
\Rightarrow&&2\left(1-\omega^2(z)\right)\omega'(z)+4\omega^2(z)\equiv 1-\omega^4(z)\\
\Rightarrow&&\omega'(z)\equiv \frac{1-4\omega^2(z)-\omega^4(z)}{2-2\omega^2(z)}\\
\Rightarrow&&4T\left(r,\omega(z)\right)=T\left(r,\omega'(z)\right)\leq 2T\left(r,\omega(z)\right)+S\left(r,\omega(z)\right)\\
\Rightarrow&&T\left(r,\omega(z)\right)=S\left(r,\omega(z)\right),\eeas
which is not possible. This completes the proof.
\end{proof}
\begin{proof}[\bf{Proof of Theorem \ref{th5}}] 
Rewritting (\ref{el2}), we get
\bea\label{el4} \left(\frac{f\left(z+c_1\right)}{f(z)}\right)^2+\left(\frac{f\left(z+c_2\right)}{f(z)}\right)^2=1.\eea 
Clearly both of $\frac{f\left(z+c_1\right)}{f(z)}$ and $\frac{f\left(z+c_2\right)}{f(z)}$ are either entire functions or non-entire meromorphic functions. We now consider the following cases.\\
{\bf Case 1.} Let $\frac{f\left(z+c_1\right)}{f(z)}$ and $\frac{f\left(z+c_2\right)}{f(z)}$ be both entire functions with $\rho_2(f)<1$. If 
$f(z)$ is any transcendental meromorphic function, then in view of \textrm{Proposition A (i)}, we claim that
\bea\label{el7a} \frac{f\left(z+c_1\right)}{f(z)}=\cos(\eta(z))\;\;
\text{and}\;\;\frac{f\left(z+c_2\right)}{f(z)}=\sin(\eta(z))
,\eea
where $\eta(z)$ is an entire function. Clearly $\eta(z)\not\equiv 0$, otherwise $f(z)$ reduces to a constant and this contradicts that $f$ is transcendental. 
If possible, let $\eta(z)$ be a non-constant entire function. If $z_0$ is a zero of $f(z)$ and it is not the zero of $f\left(z+c_1\right)$ or $f\left(z+c_2\right)$, then comparing both sides of (\ref{el7a}) we get a contradiction. So all the zeros of $f(z)$ are the zeros of $f\left(z+c_1\right)$ as well as $f\left(z+c_2\right)$. Note that $N\left(r,0;\cos(\eta(z))\right)\not=0$ and $N\left(r,0;\sin(\eta(z))\right)\not=0$. Clearly from (\ref{el7a}) we see that $N\left(r,0;f\left(z+c_1\right)\right)> N(r,0;f(z))$ as well as $N\left(r,0;f\left(z+c_2\right)\right)> N(r,0;f(z))$ and in view of \textrm{Lemma \ref{lem8l}}, we arrive at a contradiction.
So $\eta(z)$ is a non-zero constant.
Thus let $\frac{f\left(z+c_1\right)}{f(z)}\equiv M$ and $\frac{f\left(z+c_2\right)}{f(z)}\equiv N$, where $M^2+N^2=1$ with $M,N\in\mathbb{C}\setminus\{0\}$ are distinct. So 
$f(z)=k_1e^{\frac{z}{c_1}\log M}\mathscr{F}(z)=k_2e^{\frac{z}{c_2}\log N}\mathscr{G}(z)$, where $\mathscr{F}$, $\mathscr{G}$ are periodic meromorphic functions with periods $c_1$, $c_2$ respectively and $k_1,k_2\in\mathbb{C}\setminus\{0\}$,
and satisfying 
\beas \mathscr{F}\left(z+c_2\right)=\frac{k_2}{k_1}\frac{N}{M^{\frac{c_2}{c_1}}}e^{\left(\frac{\log N}{c_2}-\frac{\log M}{c_1}\right)z}\mathscr{G}(z)\;\text{and}\; 
\mathscr{G}\left(z+c_1\right)=\frac{k_1}{k_2}\frac{M}{N^{\frac{c_1}{c_2}}}e^{\left(\frac{\log M}{c_1}-\frac{\log N}{c_2}\right)z}\mathscr{F}(z).\eeas
Note that, $\rho_2\left(\mathscr{F}\right)<1$ and $\rho_2\left(\mathscr{G}\right)<1$.\\
{\bf Case 2.} Let $\frac{f\left(z+c_1\right)}{f(z)}$ and $\frac{f\left(z+c_2\right)}{f(z)}$ be both entire functions with $\rho_2(f)\geq 1$. 
Let $f(z)=e^{\phi(z)}\psi(z)$ such that $e^{\phi(z)}$ and $\psi(z)$ are respectively the non-zero part and all-zero part of $f(z)$, where $\phi(z)$ and $\psi(z)$ are entire function and meromorphic function respectively with at least one of $e^{\phi(z)}$ and $\psi(z)$ is of hyper-order $\geq 1$. We consider the following cases.\\
{\bf Sub-case 2.1.} If $\psi(z)\equiv k_1\psi\left(z+c_1\right)\equiv k_2\psi\left(z+c_2\right)$ ($k_1,k_2\in\mathbb{C}\setminus\{0\}$), then (\ref{el2}) reduces to
\bea\label{e23} k_1^{-2}e^{2\phi\left(z+c_1\right)}+k_2^{-2} e^{2\phi\left(z+c_2\right)}\equiv e^{2\phi(z)}.\eea
Suppose that $\phi(z)$ is both $c_1$  and $c_2$ periodic such that $\frac{c_1}{c_2}\in\mathbb{R}$, then we must have $ k_1^{-2}+ k_2^{-2}=1$. \\
Next suppose that $\phi(z)$ is either one of $c_1$ or $c_2$ periodic. Without loss of generality, we assume that $\phi(z)$ is $c_1$ periodic, then from (\ref{e23}), we have 
$\left(k_1^{-2}-1\right)e^{2\phi(z)}+k_2^{-2} e^{2\phi\left(z+c_2\right)}\equiv 0$. Then we must have $k_1^{-2}\not=1$ and $\phi\left(z+c_2\right)-\phi(z)$ is constant. For the rest portion, we follow \textrm{Sub-case 2.1.1} and \textrm{Sub-case 2.1.2}.\\  
Now we consider that $\phi(z)$ is neither $c_1$ nor $c_2$ periodic. So in view of \textrm{Lemma \ref{lem9}}, we conclude that at least one of $\phi\left(z+c_1\right)-\phi(z)$, $\phi\left(z+c_2\right)-\phi(z)$ and $\phi\left(z+c_1\right)-\phi\left(z+c_2\right)$ is constant. Then the following situations arise.\\
{\bf Sub-case 2.1.1.} When $\phi(z)$ is a non-zero polynomial. Let $\phi(z)=\sum_{j=0}^la_jz^j$, where $a_j\in\mathbb{C}$ ($0\leq j\leq l$). Claim that $l=1$. If possible, let $l\geq 2$. Then $\deg \left(\phi\left(z+c_1\right)-\phi\left(z+c_2\right)\right)\geq 1$ and we arise a contradiction from (\ref{e23}). Then $\phi(z)$ is such that $\phi(z)=a_1z+a_0$ with $k_1^{-2}e^{2a_1c_1}+k_2^{-2}e^{2a_1c_2}=1$, where $a_0,a_1\in\mathbb{C}$. Note that, then $\psi(z)$ reduces to the non-constant meromorphic function of hyper-order $\geq 1$.\\
{\bf Sub-case 2.1.2.} When $\phi(z)$ is a non-polynomial entire function. Let $\phi\left(z+c_1\right)-\phi(z)\equiv k$, where $k(\not=0)\in\mathbb{C}$. Now we can write $\phi$ such that $\phi(z)=\Phi_1(z)+c_4z+c_3$, where $c_3,c_4\in\mathbb{C}$ and $\Phi_1(z)$ is a transcendental entire periodic function with period $c_1$. Clearly then $c_4=\frac{k}{c_1}$ and so $\phi(z)=\Phi_1(z)+\frac{k}{c_1}z+c_3$. Finally, we can obtain the similar conclusions, if $\phi\left(z+c_2\right)-\phi(z)\equiv$ constant and $\phi\left(z+c_1\right)-\phi(z+c_2)\equiv$ constant.\\
{\bf Sub-case 2.2.} If $\psi(z)\not\equiv k_1\psi\left(z+c_1\right)$, $\psi(z)\not\equiv k_2\psi\left(z+c_2\right)$ and $k_1\psi\left(z+c_1\right)\not\equiv  k_2\psi\left(z+c_2\right)$ ($k_1,k_2\in\mathbb{C}\setminus\{0\}$), then we set by 
\bea\label{rs} f(z)=\frac{\gamma(z)}{\beta(z)}e^{\phi(z)},\eea
where $e^{\phi(z)}$, $\gamma(z)$, $\beta(z)$ are respectively non-zero part, all-zero part and all-pole part of $f(z)$. 
Clearly $\phi(z)$, $\gamma(z)$ and $\beta(z)$ are entire functions. Then (\ref{el2}) reduces to
\bea\label{en1}\frac{\gamma^2\left(z+c_1\right)}{\beta^2\left(z+c_1\right)}e^{2\phi\left(z+c_1\right)}+\frac{\gamma^2\left(z+c_2\right)}{\beta^2\left(z+c_2\right)}e^{2\phi\left(z+c_2\right)}=\frac{\gamma^2(z)}{\beta^2(z)}e^{2\phi(z)}.\eea
Let $z_p$ be a pole of $f\left(z+c_1\right)$ of multiplicity $l_p$. If $z_p$ is not a pole of $f\left(z+c_2\right)$, then $z_p$ is a pole of $f\left(z+c_1\right)^2+f\left(z+c_2\right)^2$, 
i.e., of $f^2(z)$ of multiplicities $2l_p$. Since $\frac{f\left(z+c_1\right)}{f(z)}$ and $\frac{f\left(z+c_2\right)}{f(z)}$ share $\infty$ CM, the only possibility is that $z_p$ must be a pole of $f\left(z+c_2\right)$ of multiplicities $l_p$. 
Similarly, if $z_p$ is a common pole of $f\left(z+c_1\right)$ and $f\left(z+c_2\right)$ of multiplicities $l_p$ and $l_q$ respectively, then $l_p=l_q$. Thus $f(z)$, $f\left(z+c_1\right)$ and $f\left(z+c_2\right)$ share $\infty$ CM, i.e., $\beta(z)$, $\beta\left(z+c_1\right)$ and $\beta\left(z+c_2\right)$ share $0$ CM.
For simplicity, we suppose $\beta(z)=\beta\left(z+c_1\right)=\beta\left(z+c_2\right)$. Consequently (\ref{en1}) takes the form
\bea\label{en2}\gamma^2\left(z+c_1\right)e^{2\phi\left(z+c_1\right)}+\gamma^2\left(z+c_2\right)e^{2\phi\left(z+c_2\right)}-\gamma^2(z)e^{2\phi(z)}=0.\eea
Now let $z_0$ be a common zero of $f^2\left(z+c_1\right)$ and $f^2\left(z+c_2\right)$ of multiplicities $2m_1$ and $2m_2$ respectively, where $m_1,m_2\in\mathbb{N}$. So in some neighbourhood of $z_0$, Taylor's series expansions lead to respectively 
\beas &&f^2\left(z+c_1\right)=a_{m_1}(z-z_0)^{2m_1}+a_{m_1+1}(z-z_0)^{2m_1+1}+\cdots, \;\text{where}\;a_{m_1}\not=0\\
\text{and}&& f^2\left(z+c_2\right)=a_{m_2}(z-z_0)^{2m_2}+a_{m_2+1}(z-z_0)^{2m_2+1}+\cdots, \;\text{where}\;a_{m_2}\not=0.\eeas
If $m_1< m_2$, then we see from (\ref{el2}) that
\beas\label{e3.j5} f^2(z)=a_{m_1}(z-z_0)^{2m_1}+a_{m_2}(z-z_0)^{2m_2}+\cdots.\eeas
Clearly $z_0$ is also zero of $f^2\left(z\right)$ of multiplicity $2m_1$. In particular, if $m_1=m_2$, then $f^2\left(z\right)$, $f^2\left(z+c_1\right)$ and $f^2\left(z+c_2\right)$ share $0$ at $z_0$ CM.
In another words, if $m_1< m_2$, then $\gamma\left(z\right)$, $\gamma\left(z+c_1\right)$ and $\gamma\left(z+c_2\right)$ have zero at $z_0$ of multiplicities $m_1$, $m_1$ and $m_2$ respectively. Also if $m_1= m_2$, then $\gamma\left(z\right)$, $\gamma\left(z+c_1\right)$ and $\gamma\left(z+c_2\right)$ share $0$ at $z_0$ CM. It is clear from (\ref{en2}) that this equation can be reformed like exactly of the form (\ref{en2}) after cancellation of all the common zeros among $\gamma\left(z\right)$, $\gamma\left(z+c_1\right)$ and $\gamma\left(z+c_2\right)$. Thus, without loss of generality, we may suppose that $\gamma\left(z\right)$, $\gamma\left(z+c_1\right)$ and $\gamma\left(z+c_2\right)$ are mutually prime entire functions.
We now set $\alpha(z)=\gamma(z)e^{\phi(z)}$. Then we get from (\ref{el7a}) that
\beas \frac{\alpha\left(z+c_1\right)}{\alpha(z)}=\cos{\eta(z)} \;\; \text{and}\;\; \frac{\alpha\left(z+c_2\right)}{\alpha(z)}=\sin{\eta(z)}\eeas
Since $\alpha(z)$, $\alpha\left(z+c_1\right)$, $\alpha\left(z+c_2\right)$ are mutually prime, we arrive at a contradiction.\\
{\bf Sub-case 2.3.} Let $\psi(z)\equiv k_1\psi\left(z+c_1\right)$ or $\psi(z)\equiv k_2\psi\left(z+c_2\right)$ or  $k_1\psi\left(z+c_1\right)\equiv k_2\psi\left(z+c_2\right)$ ($k_1,k_2\in\mathbb{C}\setminus\{0\}$). Let $\psi(z)\equiv k_1\psi\left(z+c_1\right)$ holds. 
As \textrm{Case 1}, we see that $\frac{f\left(z+c_1\right)}{f(z)}=k_1e^{\phi\left(z+c_1\right)-\phi(z)}\equiv \cos(\eta(z))$ and considering the zeros of both sides, we conclude that $\cos(\eta)(\not=0)\in\mathbb{C}$. 
From (\ref{el2}), we get $\frac{f\left(z+c_2\right)}{f(z)}\equiv \sin(\eta)$, where $\sin(\eta)(\not=0)\in\mathbb{C}$. The rest portion follows from \textrm{Case 1}. In this case, $\rho_2\left(\mathscr{F}\right)\geq 1$ and $\rho_2\left(\mathscr{G}\right)\geq 1$.\\
{\bf Case 3.} Let $\frac{f\left(z+c_1\right)}{f(z)}$ and $\frac{f\left(z+c_2\right)}{f(z)}$ be both non-entire meromorphic functions. Note that $f(z)\not\equiv f\left(z+c_1\right)$, $f(z)\not\equiv f\left(z+c_2\right)$ and $f\left(z+c_1\right)\not\equiv f\left(z+c_2\right)$.
In view of \textrm{Proposition B (i)}, we claim from (\ref{el4}) that
\bea\label{er1} \frac{f\left(z+c_1\right)}{f(z)}=\frac{2\omega(z)}{1+\omega^2(z)}\;\;
\text{and}\;\;\frac{f\left(z+c_2\right)}{f(z)}=\frac{1-\omega^2(z)}{1+\omega^2(z)},\eea
where $\omega(z)$ is an arbitrary non-constant meromorphic function on $\mathbb{C}$. Note that $\frac{f\left(z+c_1\right)}{f(z)}$ and $\frac{f\left(z+c_2\right)}{f(z)}$ share $\infty$ CM. 
For simplicity, we assume 
\bea\label{el7}&&\frac{f\left(z+c_1\right)}{f(z)}=\mathscr{M}(z)\\
\label{el8}\Rightarrow &&\frac{f\left(z+c_1+c_2\right)}{f\left(z+c_2\right)}=\mathscr{M}\left(z+c_2\right) \;\;\text{and}\;\; \frac{f\left(z+2c_1\right)}{f\left(z+c_1\right)}=\mathscr{M}\left(z+c_1\right) ,\eea
where $\mathscr{M}(z)$ 
is a non-entire meromorphic function. 
Now the following cases arise.\\
{\bf Sub-case 3.1.} When $\mathscr{M}(z)$ is not a periodic function of $c_1$ and $c_2$. Let $f(z)$ be of the equation (\ref{rs}) and $\alpha(z)$, $\gamma\left(z\right)$, $\gamma\left(z+c_1\right)$ and $\gamma\left(z+c_2\right)$ are as \textrm{Sub-case 2.2}. 
By $f\left(z+c_1+c_2\right)^2+f\left(z+2c_2\right)^2=f\left(z+c_2\right)^2$ or $f\left(z+2c_1\right)^2+f\left(z++c_1+c_2\right)^2=f\left(z+c_1\right)^2$, we claim that $\gamma\left(z+c_1+c_2\right)$ is co-prime entire function with both $\gamma\left(z+c_1\right)$ and $\gamma\left(z+c_2\right)$.
Let $\omega(z)=\frac{P(z)}{Q(z)}$, where $P(z)$ and $Q(z)$ are co-prime entire functions. Then we get from (\ref{er1}) that 
\begin{small}\beas &&\frac{\alpha\left(z+c_1\right)}{\alpha(z)}=\frac{2P(z)Q(z)}{P^2(z)+Q^2(z)}  \Rightarrow \frac{\alpha\left(z+c_1+c_2\right)}{\alpha\left(z+c_2\right)}=\frac{2P\left(z+c_2\right)Q\left(z+c_2\right)}{P^2\left(z+c_2\right)+ Q^2\left(z+c_2\right)} \\ 
\label{en21} \text{and} &&\frac{\alpha\left(z+c_2\right)}{\alpha(z)} = \frac{-P^2\left(z\right)+ Q^2\left(z\right)}{P^2\left(z\right)+ Q^2\left(z\right)}  \Rightarrow  \frac{\alpha\left(z+c_1+c_2\right)}{\alpha\left(z+c_1\right)} = \frac{-P^2\left(z+c_1\right)+ Q^2\left(z+c_1\right)}{P^2\left(z+c_1\right)+ Q^2\left(z+c_1\right)}.\eeas\end{small}
In view of $\alpha(z)$ and $\alpha\left(z+c_1+c_2\right)$ are co-primes with both of $\alpha\left(z+c_1\right)$ and $\alpha\left(z+c_2\right)$, we set  
\beas &&\alpha(z)\equiv e^{\xi(z)}\left(P^2(z)+Q^2(z)\right), \\
&&\alpha\left(z+c_1\right)\equiv 2e^{\xi(z)}P(z)Q(z)\equiv e^{\xi\left(z+c_1\right)}\left(P^2\left(z+c_1\right)+ Q^2\left(z+c_1\right)\right),\\
&&\alpha\left(z+c_2\right)\equiv e^{\xi(z)}\left(-P^2\left(z\right)+ Q^2\left(z\right)\right)\equiv e^{\xi\left(z+c_2\right)}\left(P^2\left(z+c_2\right)+ Q^2\left(z+c_2\right)\right),\\\text{and}
&&\alpha\left(z+c_1+c_2\right)\equiv 2P\left(z+c_2\right)Q\left(z+c_2\right)e^{\xi\left(z+c_2\right)}\\&&\equiv e^{\xi\left(z+c_1\right)}\left(-P^2\left(z+c_1\right)+ Q^2\left(z+c_1\right)\right),\eeas 
where $\xi(z)$ is an entire function.
We put $\eta(z)=\xi(z)-\xi\left(z+c_1\right)$ and $\chi(z)=\xi\left(z+c_2\right)-\xi\left(z+c_1\right)$. 
Clearly
$P\left(z+c_1\right)=\frac{1}{\sqrt{2}}\left[P(z)e^{\frac{\eta(z)}{2}}+Q(z)e^{\frac{\eta(z)}{2}}\right]$ and \\ $Q\left(z+c_1\right)=\frac{1}{\sqrt{2}i}[P(z)e^{\frac{\eta(z)}{2}}-Q(z)e^{\frac{\eta(z)}{2}}]$.
Therefore, we have 
 \beas&& \label{en19}P\left(z+c_1\right) \equiv \frac{1}{\sqrt{2}}\left[P\left(z+c_2\right)e^{\frac{\chi(z)}{2}} -Q\left(z+c_2\right)e^{\frac{\chi(z)}{2}} \right] \equiv\\&& \frac{1}{\sqrt{2}}\left[P(z)e^{\frac{\eta(z)}{2}} +Q(z)e^{\frac{\eta(z)}{2}} \right] \nonumber\\\Rightarrow 
 &&P\left(z+c_2\right)e^{\frac{\chi(z)}{2}} -Q\left(z+c_2\right)e^{\frac{\chi(z)}{2}}\equiv P(z)e^{\frac{\eta(z)}{2}} +Q(z)e^{\frac{\eta(z)}{2}}.\eeas
Similarly $P\left(z+c_2\right)e^{\frac{\chi(z)}{2}} +Q\left(z+c_2\right)e^{\frac{\chi(z)}{2}}\equiv -iP(z)e^{\frac{\eta(z)}{2}} +iQ(z)e^{\frac{\eta(z)}{2}}$.
Now we see that
\beas && P\left(z+c_2\right)=e^{\frac{\eta(z)-\chi(z)}{2}}\left(\frac{1-i}{2}P(z)+\frac{1+i}{2}Q(z)\right)\\ \text{and} 
&& -Q\left(z+c_2\right)=e^{\frac{\eta(z)-\chi(z)}{2}}\left(\frac{1+i}{2}P(z)+\frac{1-i}{2}Q(z)\right) \\ \Rightarrow 
&&P^2\left(z+c_2\right)+Q^2\left(z+c_2\right)=2P(z)Q(z)e^{\xi(z)-\xi\left(z+c_2\right)}.\eeas
Since $P^2\left(z+c_2\right)+ Q^2\left(z+c_2\right)\equiv e^{\xi(z)-\xi\left(z+c_2\right)}\left(-P^2\left(z\right)+ Q^2\left(z\right)\right)$, we have
\beas &&2P(z)Q(z)e^{\xi(z)-\xi\left(z+c_2\right)}\equiv e^{\xi(z)-\xi\left(z+c_2\right)}\left(-P^2(z)+ Q^2(z)\right)\\\Rightarrow 
&&-P^2(z)+Q^2(z)=2P(z)Q(z)\Rightarrow P(z)=\left(-1\pm\sqrt{2}\right)Q(z),\eeas
which concludes that $P(z)$ and $Q(z)$ are not co-primes and a contradiction arises\\
{\bf Sub-case 3.2.} When $\mathscr{M}(z)$ is a periodic function of $c_1$ or $c_2$. Then again as \textrm{Sub-case 2.2}, we get from (\ref{el7}) and (\ref{el8}) that 
\beas &&\frac{\alpha\left(z+c_1\right)}{\alpha(z)} \equiv \frac{\alpha\left(z+2c_1\right)}{\alpha\left(z+c_1\right)}\;\; \text{or}\;\; \frac{\alpha\left(z+c_1\right)}{\alpha(z)} \equiv \frac{\alpha\left(z+c_1+c_2\right)}{\alpha\left(z+c_2\right)} \\ \Rightarrow && \alpha(z)\equiv e^{\xi_1(z)}\alpha\left(z+c_1\right)\; \text{or} \; \alpha(z)\equiv e^{\xi_2(z)}\alpha\left(z+c_2\right), \\ &&\text{where}\; \xi_1(z)\;\text{ and}\; \xi_2(z)\;\text{ are entire functions}\\ \Rightarrow && \frac{f\left(z+c_1\right)}{f(z)} \;\;\text{or}\;\;\ \frac{f\left(z+c_2\right)}{f(z)} \;\; \text{is entire function}. \eeas
But this is not the case. This completes the proof.
\end{proof}
{\bf Conflict of Interest:} Authors declare that they have no conflict of interest. All authors made an equal contribution to the paper. Both of them read and approved the final manuscript.\\
{\bf Funding:} The second and third authors are supported by University Grants Commission (IN) fellowship and Swami Vivekananda Merit-cum-Means Scholarship (West Bengal) respectively.\\
{\bf Availability of data and materials:} Not applicable.\\
{\bf Acknowledgments:} I would also want to thank the referee and the editing team for their suggestions.


\begin{thebibliography}{33}
\bibitem{36} I. N. Baker, On a Class of Meromorphic Functions, Proceedings of the American Mathematical Society, 17(4), 819-822(Aug., 1966).
\bibitem{2} G. Ganapathy Iyer, On certain functional equations, J. Indian Math. Soc., 3, 312-315(1939).
\bibitem{3} F. Gross, On the equation $f^n+g^n =1$, Bull. Amer. Math. Soc., 72, 86-88(1966).
\bibitem{22} F. Gross, On the functional equation $f^n+g^n =h^n$, Amer. Math. Mon., 73, 1093-1096(1966).
\bibitem{50e} R. G. Halburd, R. J. Korhonen and K. Tohge, Holomorphic curves with shift-invariant hyperplanepreimages, Trans. Amer. Math. Soc., 366(8), 4267-4298(2014).
\bibitem{7} W. K. Hayman, Meromorphic functions, Clarendon Press, Oxford, 1964.
\bibitem{40} J. Heittokangas, R. Korhonen and I. Laine, On meromorphic solutions of certain nonlinear differential equations, Bull. Austral. Math. Soc., 66(2), 331-343(2002).
 
\bibitem{8} I. Laine, Nevanlinna theory and complex differential equations, De Gruyter Studies in Mathematics, 15, Walter de Gruyter \&\ Co., Berlin, 1993.
\bibitem{41} P. Li and C. C. Yang, On the nonexistence of entire solutions of certain type of nonlinear differential equations, J. Math. Anal. Appl., 320(2), 827-835(2006).
\bibitem{9} K. Liu, Meromorphic functions sharing a set with applications to difference equations, J. Math. Anal. Appl., 359(1), 384-393(2009).
\bibitem{kl} K. Liu, I Laine and L. Yang, Complex Delay-Differential Equations, Walter de Gruyter GmbH, Berlin/Boston, Vol. 78 (2021).
\bibitem{29} K. Liu and L. Z. Yang, A note on meromorphic solutions of Fermat types equations, An. Stiint. Univ. Al. I. Cuza Lasi Mat. (N. S.), 62(2)(1), 317-325(2016).
\bibitem{14} J. F. Tang and L. W. Liao, The transcendental meromorphic solutions of a certain type of non-linear differential equations, J. Math. Anal. Appl., 334(1), 517-527(2007).
\bibitem{16} C. C. Yang, A generalization of a theorem of P. Montel on entire functions, Proc. Amer. Math. Soc., 26, 332-334(1970).
\bibitem{42} C. C. Yang, On entire solutions of a certain type of nonlinear differential equation, Bull. Austral. Math. Soc., 64(3), 377-380(2001).
\bibitem{17} C. C. Yang and P. Li, On the transcendental solutions of a certain type of nonlinear differential equations, Arch. Math. (Basel), 82(5), 442-448(2004).
\bibitem{18} C. C. Yang and H. X. Yi, Uniqueness Theory of Meromorphic Functions, Science Press, Beijing/New York, 2003.
\bibitem{34} X. Zhang and L. W. Liao, On a certain type of nonlinear differential equations admitting transcendental meromorphic solutions, Sci. China, Math., 56(10), 2025-2034(2013).
\end{thebibliography}
\end{document}